\documentclass[reqno,11pt]{amsart}
\usepackage{latexsym,amssymb,amsfonts,amsmath,amsthm}
\usepackage[usenames]{xcolor}
\usepackage{eucal}
\usepackage{mathtools}
\usepackage[colorlinks,citecolor=blue,urlcolor=blue]{hyperref}
\usepackage{mathptmx} 
\newtheorem{thm}{Theorem}[section] 

\newtheorem{lem}[thm]{Lemma}

\theoremstyle{definition}

\theoremstyle{remark}
\newtheorem{rem}[thm]{Remark}
\numberwithin{equation}{section}
\newcommand{\Real}{\mathbb R}
\newcommand{\eps}{\varepsilon}

\newcommand{\F}{\mathcal{F}}

\newcommand{\E}{\mathbb{E}}

\DeclareMathOperator*{\sign}{sign}

\DeclareMathOperator{\Res}{Res}

\renewcommand{\Re}{\mathrm{Re}}
\renewcommand{\Im}{\mathrm{Im}}
\def\XXint#1#2#3{{\setbox0=\hbox{$#1{#2#3}{\int}$}
     \vcenter{\hbox{$#2#3$}}\kern-.5\wd0}}

\makeatletter
\newcommand{\doublewidetilde}[1]{{%
  \mathpalette\double@widetilde{#1}%
}}
\newcommand{\double@widetilde}[2]{%
  \sbox\z@{$\m@th#1\widetilde{#2}$}%
  \ht\z@=.9\ht\z@
  \widetilde{\box\z@}%
}
\makeatother

\begin{document}

\title[Estimation of a fractional signal in a small white noise]
{Asymptotic accuracy in estimation of a fractional signal in a small white noise}

\author{M. Kleptsyna}%
\address{Laboratoire de Statistique et Processus,
Universite du Maine,
France}
\email{marina.kleptsyna@univ-lemans.fr}

\author{D. Marushkevych}%
\address{Laboratoire de Statistique et Processus,
Universite du Maine,
France}
\email{dmytro.marushkevych.etu@univ-lemans.fr}

\author{P. Chigansky}%
\address{Department of Statistics,
The Hebrew University,
Mount Scopus, Jerusalem 91905,
Israel}
\email{pchiga@mscc.huji.ac.il}

\thanks{P. Chigansky's research was funded by ISF  1383/18 grant}
%\subjclass{}%

\dedicatory{Dedicated to the memory of Professor Robert Liptser}%
\date{\today}%

%\commby{}%
% ----------------------------------------------------------------
\begin{abstract}
This paper revisits the problem of estimating the fractional Ornstein - Uhlenbeck process observed in a linear channel with white noise of small intensity. We drive the exact asymptotic formulas for the mean square errors of the filtering and interpolation estimators. 
The asymptotic analysis is based on approximations of the eigenvalues and eigenfunctions of the signal's covariance operator.
\end{abstract}

\maketitle

\section{Introduction}
Consider the system of stochastic linear equations
\begin{equation}\label{XYsys}
\begin{aligned}
X_t =\, &  \beta \int_0^t X_s ds + B^H_t, \\
Y_t =\, & \mu \int_0^t X_s ds + \sqrt{\eps}  B_t, 
\end{aligned}
\end{equation}
where $B=(B_t,t\in \Real_+)$ and $B^H=(B^H_t, t\in \Real_+)$ are independent standard and fractional Brownian motions, 
$\beta$ and $\mu$ are constant coefficients, and $\eps$ is a small positive parameter. 

Recall that the fractional Brownian motion (fBm) with the Hurst exponent $H\in (0,1)$ is the Gaussian process with zero mean and covariance function 
$$
\E B^H_t  B^H_s  = \frac 1 2\Big(t^{2H}+s^{2H}-|t-s|^{2H}\Big), \quad s,t\in \Real_+.
$$
For $H=1/2$ the process $B^H$ coincides with the standard Brownian motion, but for all other values of $H\in (0,1)$ the fBm is 
neither a semimartingale nor a Markov process. For $H>1/2$ its increments are positively correlated and have the long range dependence 
$$
\sum_{n=1}^\infty B^H_1 (B^H_{n+1}-B^H_n)=\infty.
$$
Due to diversity of properties fBm plays an important role in the theory and applications of 
stochastic processes, see, e.g.,  \cite{M08}, \cite{PT17}.
The fractional Ornstein-Uhlenbeck (fOU) process defined by equation \eqref{XYsys} inherits the long range dependence  
from the fBm, see \cite{CKM03}.

The optimal estimation problem of the signal $X$ given the observed trajectory of the process 
$Y$ consists of computing the conditional expectation $\widehat{X}_{t,T} = \E (X_t|\F^Y_T)$ at a time  
$t\in [0,T]$ where $\F^Y_T=\sigma\{Y_t, t\le T\}$. This estimator minimizes the mean squared error over all 
functionals measurable with respect to $\F^Y_T$. For $t<T$ it is called the interpolating (or smoothing) estimator and 
for $t=T$, the filtering estimator. 

Since the process $(X,Y)$ is Gaussian, the optimal estimator is a linear functional of the observations; more precisely, it is 
the stochastic integral, see \cite{LS1},
$$
\widehat X_{t,T} = \frac 1 \mu \int_0^t h_T(s,t) dY_s,
$$
where the weight function $h_T(s,t)$ is the unique solution to the integral equation 
\begin{equation}\label{WHeq}
\eps h_T(s,t) + \int_0^T \mu^2 K(r,s) h_T(r,t)dr = \mu^2 K(s,t), \quad 0\le s\le t \le T,
\end{equation}
where the covariance kernel $K(s,t)=\E X_s X_t$ for the signal in \eqref{XYsys} has the form
\begin{equation}\label{OUKst}
K(s,t) = %\E \int_0^1 \one{v\le t }e^{-\beta(t-v) }dB^H_v \int_0^1 \one{u\le s}e^{-\beta(s-u) }dB^H_u =
\int_0^t    e^{\beta(t-v)} \frac d{dv} \int_0^s H |v-u|^{2H-1} \sign(v-u)  e^{\beta(s-u)}du dv.
\end{equation}
The minimal mean squared error $P_{t,T} (\eps)  = \E (X_t-\widehat {X}_{t,T})^2$ is given by the formula  
$$
P_{t,T} (\eps)   =
% K(t,t) - \int_0^T h(r,t) K(r,t)dr = 
\frac \eps {\mu^2} h_T(t,t).
$$ 

A closed form solution to \eqref{WHeq} is known only for $H=1/2$, i.e., when the signal is the classical OU process.
This fact lies in the foundations of the Kalman-Bucy theory of the optimal linear filtering, \cite{KB61}. 
The optimal estimate in this case can be computed recursively by solving a stochastic differential equation. The corresponding 
minimal estimation error solves the Riccati ordinary differential equation (see, e.g., \cite[Theorem 12.10 ]{LS2}), and 
a simple calculation reveals its exact asymptotics 
$$
 P_{t,T}(\eps)\asymp \sqrt{\eps/\mu^2}  \begin{cases}
1/2,  & t\in (0,T),\\
 1, & t=T,
\end{cases}\qquad \text{as\ } \eps\to 0.
$$
Here and below $f(\eps)\asymp g(\eps)$ means that $f(\eps)/g(\eps)\to 1$ as $\eps\to 0$.
This formula shows that the interpolation and the filtering errors differ asymptotically only by a constant factor and decrease as the square root of 
the noise intensity.  

\section{The main results}
In this paper we will derive a more general mean square error asymptotics, which remains valid for all values of the Hurst parameter.

\begin{thm} \label{main-thm}
The minimal mean squared estimation error of the signal in the system \eqref{XYsys} satisfies the following asymptotics as $\eps\to 0$:
\begin{equation}\label{Peps}
P_{t,T}(\eps) \asymp (\eps/\mu^2)^{\frac{2H}{1+2H}}  
 \frac{\Big(\sin (\pi H) \Gamma(2H+1)\Big)^{\frac 1 {1+2H}}}{\sin \frac \pi {2H+1}}
 \begin{cases}
 \frac{1}{2H+1}, & t\in (0,T),\\
 1, & t=T.
 \end{cases} 
\end{equation}
\end{thm}

\begin{rem}
The estimation rate $\eps^{2H/(1+2H)}$ in \eqref{Peps} coincides with the optimal minimax rate in the estimation problem of deterministic 
signals with H\"older exponent $H\in (0,1)$ in white noise, see. \cite{IKh}, \cite{Tsybakov}.
\end{rem}

\medskip

Derivation of asymptotics \eqref{Peps} uses approximations to solutions of the spectral problem  
$$
K \varphi = \lambda \varphi
$$
for the covariance operator  
$$
(K\varphi)(t) = \int_0^T K(s,t) \varphi_s ds  
$$
of the fOU process with kernel \eqref{OUKst}.
It is well known, that this problem has countably many nontrivial solutions $(\lambda_n, \varphi_n)$, $n\in \mathbb{N}$. 
The eigenvalues $\lambda_n$ are real, positive and can be ordered as a sequence decreasing to zero.  
The eigenfunctions $\varphi_n$ form an orthonormal basis in the space  
$L^2([0,T])$. 
Adding $\beta$ to the notation for the kernel \eqref{OUKst}, note that it satisfies the scaling property  
\begin{equation}\label{Kbeta}
K_\beta(sT,tT) = T^{2H} K_{\beta T}(s,t), \quad s,t\in [0,1],\quad T>0,
\end{equation}
by which the spectral problem can be considered on the unit interval $[0,1]$ without loss of generality.  

For the classical OU process, i.e. for $H=1/2$, the spectral problem can be reduced to a boundary problem for a 
linear differential equation which has simple closed form solutions
\begin{equation}
\label{OUeig}
\lambda_n = \frac 1{\nu_n^2+\beta^2}\quad  \text{and}\quad \varphi_n(t) \propto  \sqrt{2} \sin (\nu_n t),
\end{equation}
where $\nu_n = \pi n -   \pi/ 2 + O(n^{-1})$ is the increasing sequence of the roots of equation 
$$
\nu/\beta = \tan \nu.
$$
For all other values of the exponent $H\in (0,1)$ asymptotically exact approximations for the eigenvalues 
and eigenfunctions  can be obtained using the approach from \cite{Ukai}, which was recently applied in  
\cite{ChK} to the spectral analysis of the fBm. 
\medskip

The following result is of an independent interest and can be useful in other applications.  

\goodbreak

\begin{thm}\label{thm-fOU}\
\medskip

\noindent
1. The eigenvalues of the covariance operator with kernel \eqref{OUKst} on the unit interval $[0,1]$ satisfy the representation  
\begin{equation}\label{lambda_n_fOU}
\lambda_n = \sin (\pi H) \Gamma(2H+1)\frac{\nu_n^{1-2H}}{\nu_n^2+\beta^2}, \quad n=1,2,...,
\end{equation}
where the sequence $\nu_n$ has the same asymptotics as in the case of the fBm (\cite[Theorem 2.1]{ChK}):
$$
\nu_n  = \Big(n-\frac 1 2\Big)\pi -   \frac{(H-\frac 12)^2 }{H+\frac 12}\frac \pi 2+ O(n^{-1}),\quad n\to \infty.
$$

\noindent
2. The eigenfunctions with the unit norm admit the same asymptotic approximation as in the case of the fBm (\cite[Theorem 2.1 ]{ChK}):
\begin{equation}\label{fBmeigfn}
\begin{aligned}
\varphi_n& (x)  
= 
 \sqrt 2 \sin\big( \nu_{n}  x+\eta_H\big)  \\
&
 -    \int_0^{\infty}    
\Big(
e^{-  x \nu_n  u} f_0(u) +
(-1)^{n}   e^{-  (1-x) \nu_n  u}f_1(u)\
\Big)du + O(n^{-1}),
\end{aligned}
\end{equation} 
where the residual is uniform in $x\in [0,1]$, the functions $f_j(\cdot)$ are given by explicit formulas (see Lemma \ref{lem:9}) and
$$
\eta_H = \frac 1 4 \frac{(H-\frac 1 2)(H-\frac 3 2)}{H+\frac 1 2}.
$$

\medskip

\noindent 
{\bf 3.} 
The eigenfunctions satisfy    
$$
\varphi_n(1) = -(-1)^n \sqrt{2H+1} \big(1+O(n^{-1})\big), \quad n\to\infty.
$$
\end{thm}

\section{Proof of Theorem \ref{main-thm}}

Using the property \eqref{Kbeta}, we can rewrite the equation \eqref{WHeq} as
$$
\eps h(u,v) +  \int_0^1 \mu^2 T^{2H+1} K_{\beta T}(r,u)     h(r,v)dr  = 
 \mu^2 T^{2H} K_{\beta T}(u,v), \quad 0\le u\le v \le 1,
$$
where $h(u,v):= h_T(uT,vT)$. Clearly   
$$
P_{uT,T}(\eps) = \frac \eps {\mu^2} h(u,u), \quad u\in [0,1].
$$
Expanding the solution into series of the eigenfunctions of the kernel $K_{\beta T}$ gives
$$
h(u,v) = \sum_{n=1}^\infty \frac{\mu^2T^{2H}}{\eps \lambda_n^{-1} + \mu^2T^{2H+1} }
\varphi_n(u) \varphi_n(v), \quad 0\le u\le v\le 1,
$$
where $\lambda_n$ are its eigenvalues. This series is absolutely convergent for any $\eps>0$,
and its value increases unboundedly as $\eps\to 0$. Its first order asymptotics will not change if 
the eigenvalues and the eigenfunctions are replaced with their approximations from Theorem \ref{thm-fOU}. 
Denote $C:= \sin (\pi H) \Gamma(2H+1)$, then  
\begin{align*}
&
P_{T,T}(\eps) =\, 
\frac \eps {\mu^2} h(1,1)   =   
\sum_{n=1}^\infty \frac{\eps T^{2H}}{\eps \lambda_n^{-1} + \mu^2T^{2H+1} }\varphi_n^2(1) \asymp \\
&
\sum_{n=1}^\infty \frac{\eps T^{2H}}{\dfrac \eps  C 
 \dfrac{(\pi n)^2+(\beta T)^2}{(\pi n)^{1-2H}}
+ \mu^2T^{2H+1} } (2H+1) \asymp \\
&
(2H+1)\int_1^\infty \frac{\eps T^{2H}}{\dfrac \eps  C 
\big( (\pi x)^{2H+1}+(\beta T)^2 (\pi x)^{2H-1} \big)
+ \mu^2T^{2H+1} } dx \asymp \\
&
 \frac {\eps} {\mu^2T}  
\left(\frac \eps C \frac 1{T^{2H+1}\mu^2}\right)^{-\frac 1{2H+1}}
\frac{2H+1} \pi\int_0^\infty \frac{1 }{y^{2H+1}  + 1   } d y = \\
&
 \frac {\eps} {\mu^2T}  
\left(\frac \eps C \frac 1{T^{2H+1}\mu^2}\right)^{-\frac 1{2H+1}}
\frac 1 {\sin \frac \pi {2H+1}} =
 (\eps/\mu^2)^{\frac {2H}{2H+1}}    
\frac {C^{\frac 1{2H+1}}} {\sin \frac \pi {2H+1}},
\end{align*}
which is the asymptotic expression for the filtering error in \eqref{Peps} with $t=T$. 
To compute $P_{uT,T}(\eps)$ for  $u\in (0,1)$ we will use the approximation for the eigefunctions
 \eqref{fBmeigfn} which, for any fixed $u$ in the interior of the interval, takes the form  
$$
\varphi_n(u) =  \sqrt{2}\sin  \big(\nu_n u+\eta_H\big) + O(n^{-1}) \quad \text{as\ } n\to\infty.
$$ 
Then 
\begin{equation}\label{PtTu}
\begin{aligned}
&
P_{uT,T}(\eps)  =  \frac \eps {\mu^2} h(u,u) =
\sum_{n=1}^\infty \frac{ \eps T^{2H}}{\eps \lambda_n^{-1} + \mu^2T^{2H+1} }
\varphi_n^2(u)= \\
&
\sum_{n=1}^\infty \frac{\eps T^{2H}}{\eps \lambda_n^{-1} + \mu^2T^{2H+1} }
- 
\sum_{n=1}^\infty \frac{ \eps T^{2H}}{\eps \lambda_n^{-1} + \mu^2T^{2H+1} }
\cos    \big(2\nu_n u+2\eta_H\big) := \\
& I_1(\eps)+I_2(\eps).
\end{aligned}
\end{equation}
Here $I_1(\eps)$ differs from the previous case only by the factor $2H+1$, and thus, to derive the asymptotics of the interpolation error
 \eqref{Peps} for any $t\in (0,T)$, it remains to show that $I_2(\eps)$ vanishes  
as $\eps\to 0$ faster than $\eps^{2H/(2H+1)}$.
To this end, define the sequence of partial sums $(S_n)_{n\ge 0}$ 
$$
S_n = \sum_{k=1}^n \cos    \big(2\nu_k u+2\eta_H\big), \quad n\ge 1,\quad S_0=0.
$$
For any $u\in (0,1)$ this sequence is bounded and 
\begin{align*}
I_2(\eps)
= 
&
\sum_{n=1}^\infty \frac{ \eps T^{2H}}{\eps \lambda_n^{-1} + \mu^2T^{2H+1} }
\big(S_n-S_{n-1}\big) = \\
&
\eps^2 T^{2H}\sum_{n=1}^\infty S_n 
\frac
{
    \lambda_{n+1}^{-1} -\lambda_n^{-1} 
}
{
\big(\eps \lambda_{n+1}^{-1} + \mu^2T^{2H+1}\big)
\big(\eps \lambda_n^{-1} + \mu^2T^{2H+1}\big)
}.
\end{align*}
By formula \eqref{lambda_n_fOU},  for all  $n$ large enough,  we have
$$
|\lambda_{n+1}^{-1} -\lambda_n^{-1}| \le C_1 n^{2H}  \quad \text{and}\quad \lambda_n^{-1} \ge C_2 n^{2H+1}
$$
with some constants $C_1$ and $C_2$. Then, since $(S_n)$ is bounded, there exists a constant $C_3$ such that
\begin{align*}
|I_2(\eps)| \le \, 
&
C_3
\eps^2  \sum_{n=1}^\infty  
\frac
{
    n^{2H}
}
{
\big(\eps n^{2H+1} + 1\big)^2
} \asymp 
C_3
\eps^2  \int_1^\infty  
\frac
{
    x^{2H}
}
{
\big(\eps x^{2H+1} + 1\big)^2
}dx \asymp\\
&
C_3\eps \int_0^\infty  
\frac
{
   y^{2H}
}
{
\big(y^{2H+1} + 1\big)^2
}dy = O(\eps), \quad \eps\to 0.
\end{align*}
Consequently, the second term in \eqref{PtTu} is asymptotically negligible, which completes the proof.

\section{Proof of Theorem \ref{thm-fOU}}
The main idea of the proof is to reduce the spectral problem to solving a certain auxiliary system of integral 
and algebraic equations, which turns out to be more tractable for asymptotic analysis.  
A detailed description of the method appears in \cite[Section 4]{ChK}; frequently used notations and results from 
the complex analysis can be found in  monograph \cite{Gahov}. 

\subsection{The case $H>\frac 1 2$}
In this case the expression \eqref{OUKst} can be simplified by interchanging integration and derivative:
$$
K(s,t) =  \int_0^t \int_0^s e^{\beta(t-v)}e^{\beta(s-u)}c_\alpha |u-v|^{-\alpha}dudv,
$$
where we defined the new parameter $\alpha:=2-2H\in (0,1)$ and the constant $c_\alpha = (1-\frac \alpha 2)(1-\alpha)$. 
In these notations the spectral problem takes the form
\begin{equation}\label{OUeigHlarge}
\int_0^1 \left(\int_0^x \int_0^y e^{\beta(x-u)}e^{\beta(y-v)}c_\alpha |u-v|^{-\alpha}dvdu\right) \varphi(y)dy=\lambda \varphi(x),\quad x\in [0,1].
\end{equation}

\subsubsection{The Laplace transform} 
Consider the Laplace transform of the solution to equation \eqref{OUeigHlarge}
\begin{equation}\label{Laptr}
\widehat \varphi(z) = \int_0^1 e^{-zx}\varphi(x)dx, \quad z\in \mathbb{C}.
\end{equation}
Since the integral is computed over a bounded interval, it defines an entire function. 
In the following lemma, using the particular form of the kernel, we will derive an expression for 
$\widehat \varphi(z)$, central to our approach.

\begin{lem}\label{lem5.1}
Let $(\lambda, \varphi)$ solve the spectral problem \eqref{OUeigHlarge}, then the Laplace transform \eqref{Laptr}  
admits the representation  
\begin{equation}
\label{phiz} 
\widehat\varphi(z) = \widehat \varphi(-\beta)  -\frac {z+\beta}{\Lambda(z)} 
 \Big(
 \Phi_0(z)+e^{-z}\Phi_1(-z)
\Big),
\end{equation}
where 
% $\psi(0) = \int_0^1 e^{\beta x}\varphi(x)dx= \widehat \varphi(-\beta)$, 
\begin{equation}
\label{Lambda}
\Lambda(z) = \frac {\Gamma(\alpha)\lambda}{c_\alpha}(z^2-\beta^2) + \int_0^\infty \frac{2 t^{\alpha }}{t^2-z^2}dt,
\end{equation}
and the functions $\Phi_0(z)$ and $\Phi_1(z)$, defined in \eqref{Phi0Phi1}, 
are sectionally holomorphic on the cut plane $\mathbb{C}\setminus \Real_{+}$.
\end{lem}

\begin{proof} 
Differentiating both sides of equation 
\eqref{OUeigHlarge} we get
\begin{equation}\label{eigpr_fOU}
\int_0^1 \left( \int_0^{y}   e^{ -\beta v} c_\alpha |  x- v|^{-\alpha} dv\right) e^{\beta   y } \varphi(y)dy
+
\beta\lambda \varphi(x) = \lambda \varphi'(x), \quad x\in [0,1].
\end{equation}
Define the function 
\begin{equation}\label{psidef}
\psi(x) := e^{-\beta x}\displaystyle \int_x^1 e^{\beta r}\varphi(r)dr.
\end{equation}
Integration by parts gives
$$
\int_0^1 \left( \int_0^{y}   e^{ -\beta v} c_\alpha |  x- v|^{-\alpha} dv\right) e^{\beta   y } \varphi(y)dy=
\int_0^1      c_\alpha |  x- y|^{-\alpha}    \psi(y) dy.
$$
Equation \eqref{eigpr_fOU} is equivalent to the generalized spectral problem     
\begin{equation}\label{geneig}
\begin{aligned}
\int_0^1   &   c_\alpha |  x- y|^{-\alpha}    \psi(y) dy
 = \lambda \Big(\beta^2 \psi(x)-\psi''(x) \Big) 
 , \quad x\in [0,1], \\
&
\psi(1) =  0, \; 
\psi'(0)    +\beta \psi(0)  =0,
\end{aligned}
\end{equation}
where the boundary conditions follow from   definition \eqref{psidef} since 
\begin{equation}
\label{psiphi}
\psi'(x)    +\beta \psi(x)   =- \varphi(x).
\end{equation}
Plugging the identity 
$$
|x-y|^{-\alpha} =\frac 1{\Gamma(\alpha)} \int_0^\infty t^{\alpha-1}e^{-t|x-y|}dt,\quad \alpha\in (0,1),
$$
into \eqref{geneig} gives
\begin{equation}
\label{one}
  \int_0^\infty t^{\alpha-1}   u(x,t) dt
 =\frac {\Gamma(\alpha)\lambda}{c_\alpha}   \Big(\beta^2 \psi(x)-\psi''(x) \Big),
\end{equation}
where we defined the function
\begin{equation}
\label{udef}
u(x,t) := \int_0^1 \psi(y)  e^{-t|x-y|}dy. 
\end{equation}

On the other hand, integrating twice by parts and using the boundary conditions of the problem \eqref{geneig}, we get
\begin{equation}\label{psitagtag}
\begin{aligned}
\widehat \psi''(z) = & \int_0^1 \psi''(x) e^{-zx}dx = 
%\psi'(1)e^{-z}-\psi'(0) +   z\psi(1)e^{-z}-z\psi(0) + z^2 \widehat \psi(z) = \\ & 
\psi'(1)e^{-z}+(\beta  -z)\psi(0) + z^2 \widehat \psi(z).
\end{aligned}
\end{equation}
Applying the Laplace transform to \eqref{one} and substituting \eqref{psitagtag} gives 
\begin{equation}
\label{rel1}
 \int_0^\infty t^{\alpha-1}  \widehat u(z,t) dt
 =\frac {\Gamma(\alpha)\lambda}{c_\alpha}   \Big((\beta^2 -z^2)\widehat\psi(z)
 -\psi'(1)e^{-z}-(\beta  -z)\psi(0)  
  \Big).
\end{equation}
%where $\psi'(1)$ and $\psi(0)$ are related to the eigenfunction $\varphi$ in the following way:
%\begin{align*}
%& \psi(0) =  \int_0^1 e^{\beta r}\varphi(r)dr\\
%& \psi'(1)  =-  \beta \psi(1)   - \varphi(1) = -\varphi(1).
%\end{align*}

Another expression for $\widehat u(z,t)$ can be derived from definition \eqref{udef}. Differentiating it
twice we obtain the equation 
%\begin{align*}
%u'(x,t) & = % \frac d{dx}\left(\int_0^x \psi(y)  e^{-t(x-y)}dy +\int_x^1 \psi(y)  e^{-t(y-x)}dy\right) \\
%-t\int_0^x \psi(y)  e^{-t(x-y)}dy +t \int_x^1 \psi(y)  e^{-t(y-x)}dy
%\end{align*} 
\begin{equation}\label{u2tag}
u''(x,t) =  t^2  u(x,t)-2t  \psi(x)
\end{equation}
and the boundary conditions 
\begin{equation}\label{utag01}
\begin{aligned}
& u'(0,t)  = \phantom{+} t u(0,t),\\
& u'(1,t)  = -t u(1,t). 
\end{aligned}
\end{equation}
Now integrating by parts twice we get 
\begin{equation}\label{utagtag}
\begin{aligned}
\widehat u''(z,t) = & \int_0^1 u''(x,t) e^{-zx}dx = u'(1,t)e^{-z}-u'(0,t)+z\widehat u'(z,t) =\\
&
  (z-t) e^{-z} u(1,t) -(z+t) u(0,t) +z^2\widehat u(z,t),
\end{aligned}
\end{equation}
where  \eqref{utag01} was used. 
Combine the Laplace transform of \eqref{u2tag} with \eqref{utagtag} to get
\begin{equation}
\label{rel2}
\widehat u(z,t)  = \frac 1{z-t}  u(0,t)
-\frac 1{z+t}   u(1,t)e^{-z}  
- \frac{2t}{z^2-t^2} \widehat \psi(z).
\end{equation}
Plugging \eqref{rel2} into \eqref{rel1} and simplifying, we get 
$$
 \widehat \psi(z) 
 =\frac 1{\Lambda(z)} 
 \Big(
 \Phi_0(z)+e^{-z}\Phi_1(-z) 
\Big),
$$
where the function $\Lambda(z)$ is defined in   \eqref{Lambda} and 
\begin{equation}\label{Phi0Phi1}
\begin{aligned}
& \Phi_0(z):= -\frac {\Gamma(\alpha)\lambda}{c_\alpha}(\beta  -z)\psi(0)+\int_0^\infty \frac {t^{\alpha-1}}{t-z}  u(0,t)dt,\\
&
\Phi_1(z) := -\frac {\Gamma(\alpha)\lambda}{c_\alpha}\psi'(1) +  \int_0^\infty \frac {t^{\alpha-1}}{t-z}   u(1,t)  dt.
\end{aligned}
\end{equation}
Since  
$
\widehat \psi'(z) =   - \psi(0) + z \widehat \psi(z)
$
and 
$
\widehat\psi'(z)    +\beta\widehat\psi(z)   =- \widehat\varphi(z)
$
(see \eqref{psiphi}), we get  
$$
 \widehat\varphi(z)=\psi(0)- (z +\beta) \widehat\psi(z),
$$
which gives  \eqref{phiz}. 
\end{proof}

\medskip 

The next lemma elaborates the structure of the function $\Lambda(z)$ and reveal some of its useful properties.

\medskip 

\begin{lem}\label{lem4.2} \

\medskip
\noindent 
a. The function $\Lambda(z)$ admits the closed form    
\begin{equation}\label{Lfla}
\Lambda (z) =
\frac {\Gamma(\alpha)\lambda}{c_\alpha}(z^2-\beta^2)
+
 z^{\alpha-1} 
\frac{\pi   }{ \cos \frac{\pi } 2\alpha}
\begin{cases}
e^{\frac{1-\alpha}{2}\pi i}, & \arg (z) \in (0,\pi), \\
e^{-\frac{1-\alpha}{2}\pi i }, &   \arg(z)\in (-\pi, 0)
\end{cases}
\end{equation}
and has zeros at $\pm z_0 = \pm i \nu$, where $\nu>0$ is the unique real root of the equation     
\begin{equation}\label{lambdanu}
\lambda 
=
\frac {c_\alpha}{\Gamma(\alpha)} \frac {\pi   } { \cos \frac{\pi } 2\alpha}
\frac{\nu^{\alpha-1}}{\beta^2 +\nu^2}.
\end{equation}

\medskip
\noindent 
b. The limits $\Lambda^\pm(t) = \lim_{z\to t^\pm}\Lambda(z)$, where $z$ tends to $t\in \Real\setminus \{0\}$ in the 
upper and the lower half-planes, are given by the expressions 
$$
\Lambda^\pm (t) =
 \frac {\Gamma(\alpha)\lambda}{c_\alpha}(t^2-\beta^2 )
+
|t|^{\alpha-1} 
\frac{\pi   }{ \cos \frac{\pi } 2\alpha}
\begin{cases}
e^{\pm \frac{1-\alpha}{2}\pi i}, & \quad t>0 \\
e^{\mp\frac{1-\alpha}{2}\pi i }, & \quad t<0
\end{cases}
$$
and satisfy the identities 
\begin{align}
 \Lambda^+(t) & =\overline{\Lambda^-(t)}, \label{conjp}\\
 \frac{\Lambda^+(t)}{\Lambda^-(t)} & =\frac{\Lambda^-(-t)}{\Lambda^+(-t)},   \label{prop}   \\
 \big|\Lambda^+(t)\big| &  =\big|\Lambda^+(-t)\big|.   \label{absL}
\end{align}

\medskip
\noindent 
c. The argument $\theta(t):=\arg\{\Lambda^+(t)\}\in (-\pi, \pi]$ is an odd function, $\theta(-t)=-\theta(t)$,
given by the formula 
\begin{equation}
\label{thetanu}
\theta(t) = \arctan \frac{
\sin \frac{1-\alpha}{2}\pi 
}
{
\frac{(t/\nu)^2-(\beta/\nu)^2}{1+(\beta/\nu)^2}(t/\nu)^{1-\alpha} 
+
\cos \frac{1-\alpha}{2}\pi  
},\qquad t>0.
\end{equation}
It is continuous on $(0,\infty)$ with the limits  
$$
\theta(0+) := \frac{1-\alpha}{2}\pi>0\quad \text{and} \quad\theta(\infty) :=\lim_{t\to\infty}\theta(t)=0.
$$
For all sufficiently large $\nu$ the function $\theta(u;\nu):=\theta(u\nu)$ satisfies the bound     
\begin{equation}\label{thetabeta}
\Big|\theta(u;\nu) -\theta_0(u)\Big|\le g(u)(\beta/\nu)^2,
\end{equation} 
where $g(u)$ does not depend on $\nu$, is continuous on $[0,\infty)$, and grows as $g(u)\sim u^{1-\alpha}$ for $u\to 0$ and
$g(u)\sim u^{\alpha-3}$ for $u\to\infty$, and
\begin{equation}
\label{theta0}
\theta_0(u) := \lim_{\nu\to\infty}\theta(u;\nu) = 
\arctan \frac{
\sin \frac{1-\alpha}{2}\pi 
}
{
u^{3-\alpha} 
+
\cos \frac{1-\alpha}{2}\pi
}.
\end{equation}
For all $\beta\in \Real$
\begin{equation}\label{balphadef}
b_\alpha(\beta,\nu) :=\frac{1}{\pi} \int_0^\infty \theta(u; \nu) du \xrightarrow[\nu\to\infty]{} 
\frac{1}{\pi} \int_0^\infty \theta_0(u) du =
\frac{\sin{(\frac{\pi}{3-\alpha}\frac{1-\alpha}{2})}}{\sin{\frac{\pi}{3-\alpha}}}=:b_\alpha
\end{equation}
and 
\begin{equation}\label{bnu}
\big|b_\alpha(\beta,\nu)-b_\alpha\big|\le C (\beta/\nu)^2
\end{equation}
for some constant $C>0$. 

\end{lem} 

\begin{proof}\

\medskip
\noindent
a. Integration over a suitable contour shows that   
$$
\int_0^{\infty} \frac{t^\alpha}{t^2-z^2}dt  = z^{\alpha-1}\frac 1 2
\frac{\pi   }{ \cos \frac{\pi } 2\alpha}
\begin{cases}
e^{\frac{1-\alpha}{2}\pi i}, & \arg (z) \in (0,\pi) \\
e^{-\frac{1-\alpha}{2}\pi i }, &   \arg(z)\in (-\pi, 0)
\end{cases},
$$
which gives the expression \eqref{Lfla}. 
To find all zeros of $\Lambda(z)$ in the upper half-plane, let    
$z= \nu e^{i\omega}$ where $\nu>0$ and $\omega\in (0,\pi)$. Then   
equation $\Lambda(z)=0$ takes the form 
$$
\kappa (\nu^2e^{2i\omega}-\beta^2)
+ 
 \nu^{\alpha-1}e^{i(\omega-\frac \pi 2)(\alpha-1)}  =0,
$$
where we defined the constant  
$
\kappa := \displaystyle\frac {\lambda\Gamma(\alpha)}{c_\alpha}\frac{ \cos \frac{\pi } 2\alpha}{\pi   }.
$ 
The imaginary part of this equation is 
$$
\kappa  \nu^2\sin 2\omega
+ 
 \nu^{\alpha-1}\sin (\omega-\tfrac \pi 2)(\alpha-1)   =0.
$$
For all $\omega\in (0,\frac \pi 2)$ and $\omega\in (\frac \pi 2 ,\pi)$ both sines here have the same sign and hence the equality is 
possible only for $\omega = \frac \pi 2$. Thus  $\Lambda(z)$ vanishes in the upper half-plane only at the point  
$i\nu$, where $\nu$ solves the equation \eqref{lambdanu}.
By definition \eqref{Lambda}, $\Lambda(z)$ has only conjugate zeros and hence the only zero in the lower half plane is $-i\nu$. 

\medskip 
\noindent 
b. All the assertions follow directly from the explicit expression \eqref{Lfla}. 

\medskip 
\noindent
c. 
The function $f(u):= (u^2 - (\beta/\nu)^2)u^{1-\alpha}$, $u\in \Real_+$ vanishes at $u=0$ and $u=\beta/\nu$ and has the unique minimum 
$$
\min_{u\ge 0}f(u)  = -(\beta/\nu)^{3-\alpha}\frac{2}{3-\alpha}\left(\frac{1-\alpha}{3-\alpha}\right)^{\frac{1-\alpha}{2}}.
$$
Hence for a fixed $\beta$ and all sufficiently large $\nu$ the denominator  \eqref{thetanu} is bounded away from zero uniformly in $u$.   
In view of this property, all the claims are verified by a direct calculation (the value of $\beta_\alpha$ was found in \cite{ChK}).  

\end{proof}

\subsubsection{Removal of singularities} 
By (a) of Lemma \ref{lem4.2}  the expression in \eqref{phiz} has a discontinuity on the real line and two purely imaginary poles. 
Since the Laplace transform is an entire function, both of these singularities must be removable, i.e., 
the functions $\Phi_0(z)$ and $\Phi_1(z)$ must satisfy the conditions   
\begin{equation}
\label{algc}
 \Phi_0(\pm z_0)+e^{\mp z_0}\Phi_1(\mp z_0) = 0
\end{equation}
and 
$$
\lim_{z\to t^+}\frac 1 {\Lambda(z)} \big(e^{-z}\Phi_1(-z)+\Phi_0(z)\big) =
\lim_{z\to t^-}\frac 1 {\Lambda(z)} \big(e^{-z}\Phi_1(-z)+\Phi_0(z)\big), \quad t\in \Real.
$$
The latter condition can be rewritten as
\begin{align*}
&
\frac {1}{\Lambda^+(t)} 
 \Big(
 \Phi_0^+(t)+e^{-t}\Phi_1(-t)
\Big)=
\frac {1}{\Lambda^-(t)} 
 \Big(
 \Phi_0^-(t)+e^{-t}\Phi_1(-t)
\Big), \quad t>0, \\
&
\frac {1}{\Lambda^+(t)} 
 \Big(
 \Phi_0(t)+e^{-t}\Phi_1^-(-t)
\Big)=
\frac {1}{\Lambda^-(t)} 
 \Big(
 \Phi_0(t)+e^{-t}\Phi_1^+(-t)
\Big), \quad t<0
\end{align*}
or, in view of \eqref{prop}, 
\begin{equation}\label{Hp1}
\begin{aligned}
& 
\Phi_0^+(t) - \frac{\Lambda^+(t)}{\Lambda^-(t)}\Phi_0^-(t) =  e^{-t} \Phi_1(-t) \left(\frac{\Lambda^+(t)}{\Lambda^-(t)}-1\right), \\
&
\Phi_1^+(t) - \frac{\Lambda^+(t)}{\Lambda^-(t)}\Phi_1^-(t) =   e^{-t} \Phi_0(-t) \left(\frac{\Lambda^+(t)}{\Lambda^-(t)}-1\right),
\end{aligned}\quad t>0.
\end{equation}
Since $\Lambda^-(t)=\overline{\Lambda^+(t)}$ and $\theta(t)=\arg\{\Lambda^+(t)\}$, we have 
$$
\frac{\Lambda^+(t)}{\Lambda^-(t)}-1=e^{2i\theta(t)}-1 = 2i e^{i\theta(t)}\sin\theta(t).
$$
Then \eqref{Hp1} can be written as
\begin{equation}
\label{Hp}
\begin{aligned}
&
\Phi_0^+(t) - e^{2i\theta(t)}\Phi_0^-(t) = 2i  e^{-t} e^{i\theta(t)}\sin\theta(t) \Phi_1(-t),
\\
&
\Phi_1^+(t) - e^{2i\theta(t)}\Phi_1^-(t) = 2i e^{-t}  e^{i\theta(t)}\sin\theta(t) \Phi_0(-t),
\end{aligned}\qquad t>0.
\end{equation}

It follows from definition \eqref{udef} that $tu(0,t)$ and $tu(1,t)$ are bounded, and therefore, the functions defined in \eqref{Phi0Phi1} 
satisfy the a priori estimates 
\begin{equation}\label{grzero}
 \Phi_1(z)\sim z^{\alpha-1} \quad \text{and}\quad \Phi_0(z)  \sim  z^{\alpha-1} \quad \text{as}\ z\to 0,
\end{equation}
and   
\begin{equation}
\label{grinf}
\Phi_0(z)   = 2c_2(\beta-z) + O(z^{-1}) \quad \text{and}\quad
\Phi_1(z)    =  2c_1 + O(z^{-1})
 \quad \text{as\ } z\to \infty,
\end{equation}
where we defined the constants  
\begin{equation}
\label{c1c2}
c_1   = -\frac 1 2\frac{\Gamma(\alpha)\lambda}{c_\alpha}\psi'(1)  \quad \text{and}\quad 
c_2  = -\frac 1 2\frac{\Gamma(\alpha)\lambda}{c_\alpha}\psi(0).
\end{equation}

\subsubsection{Reduction to an equivalent problem }
The Laplace transform of any solution to the spectral problem \eqref{OUeigHlarge} is given by the formula  
\eqref{phiz}, where the sectionally holomorphic functions  $\Phi_0(z)$ and $\Phi_1(z)$ satisfy the estimates    
\eqref{grinf} and \eqref{grzero}, the boundary conditions \eqref{Hp} and the algebraic constraints \eqref{algc}.
Let us establish the one-to-one correspondence between all such functions and solutions to a certain system of integral 
equations on the positive real semiaxis. To this end, we will use the common technique for solving the Hilbert boundary value problem. 

Let us first consider the homogeneous Hilbert problem  of finding a function $X(z)$, sectionally holomorphic on the cut plane  
$\mathbb{C}\setminus \Real_{+}$ and satisfying the boundary conditions  
\begin{equation}
\label{Rbvphom}
X^+(t) - e^{2i\theta(t)}X^-(t) =0, \quad t\in \Real_{+}.
\end{equation}
All such functions are given by the Sokhotski-Plemelj formula:
\begin{equation}
\label{Xz}
X(z)= z^k X_c(z)=z^k \exp \left(\frac 1 \pi \int_0^\infty\frac{\theta(t)}{t-z}dt\right), \quad z\in \mathbb{C}\setminus \Real_{+},
\end{equation}
where $k$ is an integer, to be chosen later. The canonical part $X_c(z)$ of this expression satisfies the estimates  
\begin{equation}\label{Xczinf}
X_c(z) = 1 - z^{-1} \nu b_\alpha(\beta, \nu) + O(z^{-2}) \quad \text{as}\ z\to\infty,
\end{equation}
where $b_\alpha(\beta, \nu)$ is defined in \eqref{balphadef}, and
\begin{equation}
\label{Xcz0}
X_c(z) \sim z^{ \frac{\alpha-1}{2}}  \quad \text{as}\ z\to 0.
\end{equation}

Define the functions  
\begin{equation}
\label{SDdef}
\begin{aligned}
& S(z):= \frac{\Phi_0(z)+\Phi_1(z)}{2X(z)}, \\
& D(z):= \frac{\Phi_0(z)-\Phi_1(z)}{2X(z)},
\end{aligned}
\end{equation}
which, in view of \eqref{Hp} and \eqref{Rbvphom}, satisfy the {\em decoupled} boundary conditions
\begin{equation}\label{bcndSD}
\begin{aligned}
S^+(t)-S^-(t)  &= \phantom{+}2i h(t) e^{-t} S(-t), \\
D^+(t)-D^-(t)  &=  - 2i h(t)e^{-t} D(-t),
\end{aligned}\qquad t>0,
\end{equation}
where   
$$
h(t):=e^{i\theta(t)}\sin \theta(t)\frac{X(-t)}{X^+(t)}.
$$
Calculations, similar to \cite[eq. (5.37)]{ChK}, show that this function can be written as      
\begin{equation}\label{ht}
h(t) = \exp \left(-\frac 1 \pi \int_0^\infty\theta'(s) \log\left| \frac{t+s}{t-s}\right|ds\right)\sin \theta(t),
\end{equation}
and therefore, satisfies the H\"older property on  $\Real_{+}$ and has the limit  
$$h(0):=\sin \theta(0+)=\sin \frac{1-\alpha}{2} \pi.$$

Applying the Sokhotski-Plemelj formula  to \eqref{bcndSD}, we obtain the following representation for the functions in \eqref{SDdef}:
\begin{equation}
\label{SD}
\begin{aligned}
S(z) = &\phantom{+}\frac 1 \pi \int_0^\infty \frac{h(t)e^{-t}}{t-z}S(-t) dt + P_S(z), \\
D(z) = &-\frac 1 \pi \int_0^\infty \frac{h(t)e^{-t}}{t-z}D(-t)dt + P_D(z),
\end{aligned}
\end{equation}
where the polynomials $P_S(z)$ and $P_D(z)$ are chosen to match the a priori growth estimates for $S(z)$ and $D(z)$ as $z\to\infty$. 
Note that the integrals in the right hand side  of \eqref{SD} are well defined and finite, only if  $S(-t)$ and $D(-t)$ are integrable at zero.
In view of the estimates \eqref{grzero} and \eqref{Xcz0}, this limits the choice of the integer $k$ in the expression \eqref{Xz} 
to $k < (\alpha +1)/2$. In what follows we will need $S(-t)$ and $D(-t)$ to be square integrable, 
which reduces the limitation further to $k < \alpha/ 2$. A convenient choice is $k=0$, which 
corresponds to setting $X(z):=X_c(z)$. 

Since for any real numbers $a$, $b$ and $c$
$$
\frac {az+b+O(z^{-1})}{1-cz^{-1}+ O(z^{-2})}=a(z+c) + b + O(z^{-1})\quad \text{as} \ z\to\infty,
$$
the a priori estimates \eqref{grinf} and \eqref{Xczinf} give 
\begin{align*}
S(z) & =  c_2 \big(-z+\beta -\nu b_\alpha(\beta,\nu)\big) + c_1 + O(z^{-1}),\\
D(z) & =  c_2 \big(-z+\beta -\nu b_\alpha(\beta,\nu)\big) - c_1 + O(z^{-1}).
\end{align*}
This asymptotics determines the choice of the polynomials in \eqref{SD}:
\begin{align*}
P_S(z) &:=  c_2(-z+\beta-\nu b_\alpha(\beta,\nu))+c_1,\\
P_D(z) &:=  c_2(-z+\beta-\nu b_\alpha(\beta,\nu))-c_1, 
\end{align*}
where the constants $c_1$ and $c_2$ are defined by \eqref{c1c2}. 
Now, plugging $z:=-t$ with $t\in \Real_+$ into \eqref{SD}, we obtain integral equations for the restrictions $S(-t)$ and $D(-t)$:
$$
\begin{aligned}
& 
S(-t) = \phantom{+}\frac 1 \pi \int_0^\infty \frac{h(s)e^{-s}}{s+t}S(-s) ds + c_2\big(t+\beta-\nu b_\alpha(\beta,\nu)\big)+c_1,\\
&
D(-t) = -\frac 1 \pi \int_0^\infty \frac{h(s)e^{-s}}{s+t}D(-s)ds + c_2\big(t+\beta-\nu b_\alpha(\beta,\nu)\big)-c_1.
\end{aligned}
$$
Consider the auxiliary integral equations 
\begin{equation}\label{qp}
p^\pm_j (t)  = \pm \frac 1 \pi \int_0^\infty \frac{h_\beta(s;\nu)e^{-\nu s}}{s+t} p^\pm_j(s)ds+t^j, \quad j\in \{0,1\},
\end{equation}
where $h_\beta(u;\nu):=h(u \nu)$, $u>0$. Below we will show, that for all sufficiently large $\nu$ these equations have 
the unique solutions in the class of functions such that $p^\pm_j (t)-t^j$ belong to the space $L^2(\Real_+)$. 
We will extend their domain the cut plane by replacing $t$ in the right hand side of the equation   
\eqref{qp} with $z\in \mathbb{C}\setminus \Real_{+}$.

Since by construction $S(-t)$ and $D(-t)$ are square integrable at zero, by linearity   
\begin{align*}
& 
S(z\nu ) =   
  c_2 \nu p^+_1(-z) + \Big(c_2\big(\beta-\nu b_\alpha(\beta,\nu)\big)+c_1\Big)p^+_0(-z),\\
&
D(z\nu ) =  
 c_2 \nu p^-_1(-z) + \Big(c_2\big(\beta-\nu b_\alpha(\beta,\nu)\big)-c_1\Big)p^-_0(-z).
\end{align*}
Thus letting 
\begin{align*}
&
a_\pm (z):= p^+_0(z)\pm p^-_0(z), \\
&
b_\pm (z):= p^+_1(z)\pm p^-_1(z)
\end{align*} 
and using definition \eqref{SDdef} we get 
\begin{equation}\label{Phiba}
\begin{aligned}
\Phi_0(z\nu) & =
c_2 \nu X(z\nu)\Big(b_+(-z) +   \big(  \beta/\nu -  b_\alpha(\beta,\nu)\big)a_+(-z)\Big)  +c_1X(z\nu) a_-(-z),
\\
\Phi_1(z\nu) & = 
c_2 \nu X(z\nu)\Big(b_-(-z) +   \big(\beta/\nu-  b_\alpha(\beta,\nu)\big)a_-(-z)\Big)+c_1X(z\nu) a_+(-z).
\end{aligned}
\end{equation} 
Plugging this expression into condition \eqref{algc} to get 
\begin{equation}\label{linsys}
c_2 \nu \xi
 + 
c_1 \eta  =0,
\end{equation} 
where $X_\beta(z;\nu):=X(z\nu)$ and
\begin{align}
\label{xieta}
\xi :=\, & 
e^{i\nu/2}
X_\beta(i;\nu)\Big(b_+(-i) +   \big(  \beta/\nu -  b_\alpha(\beta,\nu)\big)a_+(-i)\Big) + \\
\nonumber & \hskip 2.5cm  e^{-i\nu/2}
  X_\beta(-i;\nu)\Big(b_-(i) +   \big(\beta/\nu-  b_\alpha(\beta,\nu)\big)a_-(i)\Big),
\\  
\nonumber
\eta :=\,  &
e^{i\nu/2} X_\beta(i;\nu) a_-(-i)+e^{-i\nu/2}  X_\beta(-i;\nu) a_+(i).
\end{align}
Since the constants $c_1$ and $c_2$ are real, the linear algebraic system \eqref{linsys} has nontrivial solutions if and only if  
\begin{equation}
\label{Im}
\Im\{\xi \overline{\eta}\}=0,
\end{equation}
in which case $c_1 = - c_2 \nu \xi /\eta$. Plugging this equality into \eqref{Phiba} we get 
\begin{equation}
\label{Phinorm}
\begin{aligned}
\Phi_0(z) /c_2 \nu  = &
 X(z)\Big(b_+(-z/\nu) +   \big(  \beta/\nu -  b_\alpha(\beta,\nu)\big)a_+(-z/\nu)\Big)  \\
 & -   \frac \xi \eta X(z) a_-(-z/\nu),
\\
\Phi_1(z)/c_2 \nu  = &
  X(z)\Big(b_-(-z/\nu) +   \big(\beta/\nu-  b_\alpha(\beta,\nu)\big)a_-(-z/\nu)\Big) \\
  &-   \frac \xi  \eta X(z) a_+(-z/\nu).
\end{aligned}
\end{equation}

To recap, we arrive at a problem equivalent to solving the equation \eqref{OUeigHlarge}.

\begin{lem}\label{lem5.3}
Let $(p^\pm_0, p^\pm_1, \nu)$ be a solution to the system which consists of the integral equations  \eqref{qp}  
and the algebraic conditions \eqref{Im}. Define $\varphi$ by the Laplace transform, given by the expression \eqref{phiz}, 
where $\Phi_0(z)$ and $\Phi_1(z)$ are defined by the formulas \eqref{Phinorm}, and the number $\lambda\in \Real_+$ by the formula 
\eqref{lambdanu}. Then the pair $(\lambda, \varphi)$ solves the spectral problem  \eqref{OUeigHlarge}. 
Conversely, starting with a solution $(\lambda, \varphi)$ to the problem \eqref{OUeigHlarge}, a solution to the above 
integro-algebraic system can be constructed.  
\end{lem}

The following lemma derives the exact asymptotics for  $X_\beta(i;\nu)$ as $\nu\to \infty$. 
This limit will be used in calculations to follow. 

\begin{lem}\label{lemXbeta}
$$
\begin{aligned}
\arg\big\{X_\beta(i;\nu)\big\} = & \frac{1-\alpha}{8}\pi + O(\nu^{-2}), \\
\big|X_\beta(i;\nu)\big| = & \sqrt{\frac{3-\alpha}{2}} + O(\nu^{-2}),
\end{aligned}\qquad \nu\to\infty.
$$
\end{lem}

\begin{proof}
The constants in the right  hand side are the argument and the absolute value of the limit, cf.  \eqref{theta0},
$$
X_0(i):= \lim_{\nu\to\infty}X_\beta(i;\nu) = 
\exp \left(\frac 1 \pi \int_0^\infty\frac{\theta_0(u)}{u -i}du\right),
$$ 
calculated in   \cite[Lemma 5.5]{ChK}. The estimates for the residuals follow from inequality   
\eqref{thetabeta}.
\end{proof}

\subsubsection{Properties of the integro-algebraic system}
Solvability of the system, introduced in Lemma \ref{lem5.3}, is guaranteed by contractivity of the operator    
$$
(A f)(t) := \frac 1\pi \int_0^\infty \frac{h_\beta(s;\nu)e^{-\nu s}}{s+t}f(s)ds,
$$
where $h_\beta(u;\nu):= h(u\nu)$ (see \eqref{ht}):
\begin{equation}
\label{hbetanu}
h_\beta(u;\nu) = 
\exp \left(-\frac 1 \pi  \int_0^\infty\theta'(v;\nu) \log\left| \frac{u +v }{u -v }\right|dv \right)
\sin \theta(u;\nu).
\end{equation}

\begin{lem}\label{lem5.5}
The operator $A$ is a contraction on the space $L^2(\Real_+)$ for all sufficiently large $\nu$, i.e., 
for any $\alpha_0\in (0,1]$ there exist constants $\eps>0$ and  $\nu'>0$, such that
$\|A\|\le 1-\eps$ for all $\nu\ge \nu'$ and $\alpha \in [\alpha_0,1]$. 
\end{lem}

\begin{proof}
A direct calculation shows that for all sufficiently large  $\nu$ and all $\alpha\in [\alpha_0,1]$ 
the exponent in \eqref{hbetanu} is bounded by a continuous function  $f(u)$, which does not depend on $\alpha$ 
and $\nu$, and whose limits as $u\to 0$ and $u\to\infty$ equal 1. Thus  
$$
h_\beta(u,\nu)\le f(u) \sin \theta(u;\nu)\le \|f\|_\infty.
$$
% Here is the calculation, mentioned above.
%Consider the derivative 
%$$
%\theta'(u;\nu) = 
%\frac {-\sin \frac{1-\alpha}{2}\pi} {(f(u) 
%+
%\cos \frac{1-\alpha}{2}\pi)^2+ 
%\sin^2 \frac{1-\alpha}{2}\pi
%}f'(u)
%$$
%where we defined  $f(u):= \frac{u^2-(\beta/\nu)^2}{1+(\beta/\nu)^2 }u^{1-\alpha}$, which satisfies 
%$$
%\min_{u\ge 0} f(u) = -\frac{(\beta/\nu)^{3-\alpha}}{1+(\beta/\nu)^2 }\frac{2}{3-\alpha}\left(\frac{1-\alpha}{3-\alpha}\right)^{\frac{1-\alpha}{2}} \ge - (\beta/\nu)^{3-\alpha}, \quad \alpha \in (0,1).
%$$
%Let $\nu'$ be such that $(\beta/\nu)^2 \le \frac 1 2 \cos \frac{1-\alpha_0}{2}\pi$ for all $\nu\ge \nu'$, so that
%$f(u)\ge -\frac 1 2 \cos \frac{1-\alpha_0}{2}\pi$ for all $u\ge 0$. Consequently for $u\le 1$, we have 
%$$
%\big|\theta'(u;\nu)\big| \le  
%\frac {1} {\frac 1 4  \cos^2 \frac{1-\alpha}{2}\pi + 
%\sin^2 \frac{1-\alpha}{2}\pi
%}|f'(u)|\le 4 |f'(u)| \le  24  u^{-1}
%$$ 
%Similarly, for $u>1$ we have 
%$
%f(u)
%%= \frac{u^2-(\beta/\nu)^2}{1+(\beta/\nu)^2 }u^{1-\alpha}
%\ge 
%\frac{1 }{4}u^{3-\alpha}
%$
%and therefore 
%$$
%\big|\theta'(u;\nu)\big|\le 
%16\, u^{2\alpha-6}|f'(u)| \le 
%64\, u^{ \alpha-4} \le 64\, u^{ -3}.
%$$
%By the formula \eqref{hbetanu} 
%$$
%\big|h_\beta(u;\nu)\big| \le g_0(u)
%\big|\sin \theta(u;\nu)\big| =: \widetilde h_\beta(u;\nu)
%$$
%where 
%$$
%g_0(u)= \exp \left(64\int_0^\infty (v^{-3}\wedge v^{-1}) \log\left| \frac{u+v}{u-v}\right|dv\right).
%$$
%
Also, since $\theta_0(0+) = \frac {1-\alpha}2\pi$, in view of estimate \eqref{thetabeta}, there exists a neighborhood of the origin, 
where for all sufficiently large $\nu$   
$$
\sin \theta(u;\nu) \le \frac 1 2+\frac 1 2\sin \frac {1-\alpha_0}2\pi =: 1-3\eps.
$$ 
Since $f(0)=1$, this also guarantees that $h_\beta(u,\nu)<1-2\eps$ in some neighborhood of the origin.
Thus, as $\sup_{\nu>0}\|h_\beta\|_\infty\le \|f\|_\infty$, a constant $\nu'$ can be chosen so that  
$h_\beta(u,\nu)e^{-\nu u}<1-\eps$ for all $u>0$ and all  $\nu\ge \nu'$. 
The rest of the proof can be repeated as in \cite[Lemma 5.6]{ChK}.

\end{proof}

The following estimates play the key role in the asymptotic analysis of the integro-algebraic system from Lemma \ref{lem5.3}.

\begin{lem}\label{lemest}
For any $\alpha_0\in (0,1]$ there exist constants $\nu'$ and $C$, such  that for all $\nu\ge \nu'$ and 
$\alpha\in [\alpha_0,1]$
\begin{align*}
& \big|a_-(\pm i)\big| \le C\nu^{-1}, \quad \big|a_+(\pm i)-2\big|\le C \nu^{-1}, \\
& \big|b_-(\pm i)\big| \le C\nu^{-2}, \quad \big|b_+(\pm i)\mp 2i\big|\le C\nu^{-2}, 
\end{align*}
and, for all $\tau>0$, 
\begin{align*}
& \big|a_-(\tau)\big| \le C\nu^{-1}\tau^{-1}, \quad \big|a_+(\tau)-2\big|\le C \nu^{-1}\tau^{-1}, \\
& \big|b_-(\tau)\big| \le C\nu^{-2}\tau^{-1}, \quad \big|b_+(\tau)\mp 2\tau\big|\le C\nu^{-2} \tau^{-1}.
\end{align*}
\end{lem}

\begin{proof}
As shown in the proof of the previous lemma, the function  $h_\beta(u;\nu)$ is bounded by a constant, which depends only 
on $\alpha_0$, for all $\nu$ large enough. In view of this estimate, the proof in \cite[Lemma 5.7]{ChK} applies without any changes. 
\end{proof}

\subsubsection{Laplace transform inversion}
The next lemma expresses the eigenfunctions in  terms of the solutions to integro-algebraic system from Lemma
 \ref{lem5.3}.

\begin{lem}\label{lemeigf}
Let $(\Phi_0,\Phi_1,\nu)$ be the solutions to  integro-algebraic system from Lemma \ref{lem5.3}. Then 
the function $\varphi$, defined by the Laplace transform \eqref{phiz}  satisfies  
\begin{align}\nonumber 
\varphi(x)  = - & 
\nu^{3-\alpha}\frac{ \cos \frac{\pi } 2\alpha}{\pi   }   2
\Re\bigg\{e^{i\nu x}\Phi_0(i\nu )\frac { 1-i(\beta/\nu) }{  
\frac{2}{(\beta/\nu)^2 +1}  -\alpha+1 }\bigg\}+\\
\label{phifla}
&
  \nu^{3-\alpha}\frac{ \cos \frac{\pi } 2\alpha}{\pi   }  \frac 1 \pi\int_0^\infty \frac{  \sin\theta_\beta(u;\nu)}{\gamma_\beta(u;\nu)}\cdot 
\\
&
\nonumber  
\cdot \bigg( e^{-(1-x)u\nu} \big(u+\tfrac \beta \nu\big) \Phi_1(-u\nu)
 - 
 e^{-u\nu x} \big(u-\tfrac \beta \nu \big)\Phi_0(-u\nu) \bigg)du,
\end{align}
where $\gamma_\beta(u;\nu)$ is defined by formula \eqref{gammabeta} below.
Moreover, the following equalities hold
\begin{equation}\label{flaphi}
\begin{aligned}
\int_0^1 e^{\beta x}\varphi(x)dx &= -\nu^{3-\alpha}\frac {\cos \frac \pi 2 \alpha}{\pi}  2 c_2  \big(1+(\beta/\nu)^2\big), \\
\varphi(1) &= -\nu^{3-\alpha}\frac {\cos \frac \pi 2 \alpha}{\pi} 2 c_2\nu \frac \xi\eta \big(1+(\beta/\nu)^2\big).
\end{aligned}
\end{equation}

\end{lem}

\begin{proof}
Since $\widehat \varphi(z)$ is holomorphic, the Laplace transform \eqref{phiz} inversion can be done by means of integration along the 
imaginary axis:
$$
\begin{aligned}
\varphi(x) & =  -\frac 1{2\pi i}\lim_{R\to\infty}\int_{-iR}^{iR} 
\left( 
   \frac {z+\beta}{\Lambda(z)}
 \Phi_0(z)+ \frac {z+\beta}{\Lambda(z)} e^{-z}\Phi_1(-z)
 -\psi(0)
\right)e^{zx}dz \\
& =
-\frac 1{2\pi i}\lim_{R\to\infty}\int_{-iR}^{iR} \big(f_0(z) +  f_1(z)\big)dz,
\end{aligned}
$$
where we defined the functions
$$
f_0(z) = e^{zx}\left((z+\beta)\frac {\Phi_0(z)}{\Lambda(z)}   -\psi(0)\right)\quad \text{and}\quad 
 f_1(z) = e^{(x-1)z} (z+\beta) \frac {\Phi_1(-z)}{\Lambda(z)}.
$$
Integrating on suitable contours as in the proof of \cite[Lemma 5.8]{ChK}, we get  
\begin{align*}
\int_{-i\infty}^{i\infty} \big(f_1(z)+f_0(z)\big)dz &=\,  2\pi i \Big(\Res(f_0,z_0) +  \Res(f_0,-z_0)\Big) +\\
&
\int_0^\infty \big(f_1^+(t)-f_1^-(t)\big)dt
+\int_0^\infty \big(f_0^-(-t)-f_0^+(-t)\big)dt.
\end{align*}
Due to symmetries \eqref{conjp} and \eqref{absL} and the definition of $\theta(t)$,
\begin{align*}
&
f_1^+(t)-f_1^-(t)   %= e^{(x-1)t} (t+\beta) \Phi_1(-t) \left(\frac {1}{\Lambda^+(t)} -  \frac {1}{\Lambda^-(t)}\right)
= -e^{(x-1)t} (t+\beta) \Phi_1(-t)\frac{2i \sin\theta(t)}{\gamma(t)}, \\
&
f_0^-(-t)-f_0^+(-t) %=e^{-tx} (-t+\beta)\Phi_0(-t)\left(\frac {1}{\Lambda^-(-t)}  - \frac {1}{\Lambda^+(-t)} \right)
=-e^{-tx} (-t+\beta)\Phi_0(-t)\frac{2i \sin\theta(t)}{\gamma(t)},
\end{align*}
with $\gamma(t)=|\Lambda^+(t)|$, and therefore,   
\begin{align*}
\varphi(x) = &- \Res\big(f_0,z_0\big) -  \Res\big(f_0,-z_0\big) +\\
&
 \frac 1 \pi \int_0^\infty \frac{  \sin\theta(t)}{\gamma(t)}\left( e^{-(1-x)t} (t+\beta) \Phi_1(-t)-
 e^{-tx} (t-\beta)\Phi_0(-t) \right)dt.
\end{align*}
The usual calculations produce the following expressions for the residues 
\begin{align*}
\Res\big(f_0,z_0\big)  =\, & %\lim_{z\to z_0}(z-z_0) f_0(z) =
e^{i\nu x} (i\nu +\beta)\frac {\Phi_0(i\nu )}{\Lambda'(i\nu ) }  = \\
&
e^{i\nu x}\Phi_0(i\nu ) \frac  { \cos \frac{\pi } 2\alpha}{\pi   }\nu^{3-\alpha} \frac { 1-i(\beta/\nu) }{  
\frac{2}{(\beta/\nu)^2 +1}  -\alpha+1 }, \\
\Res\big(f_0,-z_0\big)  =\, &  e^{-i\nu x} (-i\nu +\beta)\frac {\Phi_0(-i\nu )}{\Lambda'(-i\nu)}= \\
&
e^{-i\nu x}\Phi_0(-i\nu )\frac  { \cos \frac{\pi } 2\alpha}{\pi   }\nu^{3-\alpha}\frac {1  +i(\beta/\nu)}{
\frac{2 }{(\beta/\nu)^2 +1}   -\alpha+1},
\end{align*}
and consequently 
\begin{align*}
&
\Res\big(f_0,z_0\big)+\Res\big(f_0,-z_0\big) = \\
&
2\nu^{3-\alpha}\frac  { \cos \frac{\pi } 2\alpha}{\pi   }  
\Re\left\{e^{i\nu x}\Phi_0(i\nu )\frac { 1-i(\beta/\nu) }{  
\frac{2}{(\beta/\nu)^2 +1}  -\alpha+1 }\right\}. 
\end{align*}
By plugging this expression, we get \eqref{phifla}  where 
\begin{equation}
\label{gammabeta}
\gamma_\beta(u;\nu)=\nu^{1-\alpha}\frac{ \cos \frac{\pi } 2\alpha}{\pi   } \big|\Lambda^+(u\nu)\big|= \left| 
\frac{u^2-(\beta/\nu)^2 }{(\beta/\nu)^2 +1}
+
u^{\alpha-1} 
e^{ \frac{1-\alpha}{2}\pi i}\right|.
\end{equation}
The formulas \eqref{flaphi} follow  from \eqref{c1c2}, \eqref{lambdanu} and \eqref{linsys}. 
\end{proof}

\subsubsection{Asymptotic analysis}
By Lemma \ref{lem5.3} the spectral problem \eqref{OUeigHlarge} reduces to solving the integro-algebraic system of equations.
The following lemma gives the exact asymptotics of its algebraic part.

\begin{lem}\label{lem5.8}
The integro-algebraic system of Lemma  \ref{lem5.3} has countably many solutions, which can be enumerated so that 
\begin{equation}\label{nuas}
\nu_n = \pi \Big(n+\frac 1 2\Big) -\frac{1-\alpha}{4}\pi + \arcsin{\frac{b_\alpha}{\sqrt{1+b^2_\alpha}}} + n^{-1} r_n(\alpha), \quad n\to \infty,
\end{equation}
where the residual  $r_n(\alpha)$ is uniformly bounded in $n\in \mathbb{N}$ and $\alpha \in [\alpha_0, 1]$ for all $\alpha_0\in (0,1]$.
\end{lem}

\begin{proof}
The proof is similar to  \cite[Lemma 5.9]{ChK}.
Plugging the estimates from Lemmas \ref{lemest} and \ref{lemXbeta} and the bound \eqref{bnu} into definition  
\eqref{xieta}, we can write  
$$
\xi \overline{\eta} = 4 \frac{3-\alpha}{2}\sqrt{1+b_\alpha^2}\exp \left\{i\Big(\nu+ \frac{1-\alpha}{4}\pi-\pi+\arg\big\{i+b_\alpha\big\}\Big)\right\}
\big(1+R(\nu)\big),
$$
where the function $R(\nu)$ satisfies the inequality $|R(\nu)|\le C_1 \nu^{-1}$ with some constant $C_1$ which depends only on $\alpha_0$. 
Hence equation \eqref{Im} takes the form 
\begin{equation}
\label{Imexplicit}
\nu+ \frac{1-\alpha}{4}\pi-\pi+\arg\big\{i+b_\alpha\big\}-\pi n+\arctan \frac{\Im\{R(\nu\}}{1+\Re\{R(\nu)\}}=0,  
\end{equation}
for all  $n\in \mathbb{Z}$.
This fixes a certain enumeration of all the solutions to the integro-algebraic system of Lemma \ref{lem5.3}. Clearly,   $\nu$ is positive 
for all $n$ large enough. However at this point solvability of this equation for any such $n$ is not obvious. This can be argued 
as follows.

By Lemma \ref{lem5.5} the integral operator in the right hand side of equations \eqref{qp} is contracting in $L_2(\Real_+)$ 
for all sufficiently large $\nu$.  A direct calculation shows that $|R'(\nu)|\le C_2\nu^{-1}$ with some constant   
$C_2$. Hence for all $n$ large enough the system which consists of the integral  and algebraic equations, 
\eqref{qp} and \eqref{Imexplicit}, has the unique solution obtained by fixed point iterations of the integro-algebraic operator.
Asymptotics \eqref{nuas} follows from \eqref{Imexplicit}, since 
$$
\arg\{i+b_\alpha\}=\frac \pi 2 -\arcsin \frac{b_\alpha}{\sqrt{1+b_\alpha^2}}.
$$
\end{proof}

The next lemma derives the corresponding asymptotic approximation of the eigenfunctions.

\begin{lem} \label{lem:9}
The eigenfunctions, enumerated as in Lemma \ref{lem5.8},  admits the approximation 
\begin{align}\label{phias}
&
\varphi_n(x) =   \sqrt{2}
\cos \Big(\nu_n x + \frac {1-\alpha}8 \pi + 
\frac \pi 2 -\arcsin \frac{b_\alpha}{\sqrt{1+b_\alpha^2}}
\Big)
\\
\nonumber
& +
 \frac {\sqrt{3-\alpha}} \pi\int_0^\infty \rho_0(u)
 \Big(    -e^{-u\nu_n x}  \frac{u-b_\alpha}{\sqrt{1+b_\alpha^2}} -(-1)^ne^{-(1-x)u\nu_n} \Big)du + n^{-1}r_n(x),
\end{align}
where the residual $r_n(x)$ is uniformly bounded in   $n\in \mathbb{N}$ and $x\in [0,1]$, and    
$$
\rho_0(u)=\frac{  \sin\theta_0(u)}{\gamma_0(u)}X_0(-u).
$$
Moreover, 
\begin{equation}
\label{psias}
\begin{aligned}
\varphi_n(1) =  &
- (-1)^n \sqrt{3-\alpha}  \big(1+O(n^{-1})\big),
\\
\int_0^1 e^{\beta x}\varphi_n(x)dx = &  
- \sqrt{\frac{3-\alpha}{1+b_\alpha^2}} \nu_n^{-1}
\end{aligned}
\end{equation}
and 
\begin{equation}\label{Dima}
\int_0^1 \varphi_n(x)dx  =- \sqrt{\frac{3-\alpha}{1+b_\alpha^2}} \nu_n^{-1}.
\end{equation}
\end{lem}

\begin{proof} 
Let  
$
\gamma_0(u):= \big|u +u^{\alpha-2} e^{ \frac{1-\alpha}{2}\pi i}\big|,
$
then, due to \eqref{gammabeta},
$$
\big|\gamma_\beta(u;\nu)-u\gamma_0(u)\big|\le 2(\beta/\nu)^2(u^2+1).
$$
Along with \eqref{thetabeta} expression \eqref{phifla} gives
\begin{multline*}
\varphi_n(x) \propto  - \frac 2 {3-\alpha}\Re\Big\{e^{i\nu_n x}\Phi_0(i\nu_n )\Big\}\\
+
 \frac 1 \pi\int_0^\infty \frac{  \sin\theta_0(u)}{\gamma_0(u)}\left( e^{-(1-x)u\nu_n}   \Phi_1(-u\nu_n)-
 e^{-u\nu_n x}  \Phi_0(-u\nu_n) \right)du + n^{-1}r_n(x)
\end{multline*}
where the residual  $r_n(x)$ is uniformly bounded in $n\in \mathbb{N}$ and $x\in [0,1]$. Approximation \eqref{phias} is obtained 
by plugging the estimates from Lemmas  \ref{lemest} and \ref{lemXbeta}, and the bound \eqref{bnu}  into expression \eqref{Phinorm} and 
normalizing to unit $L^2([0,1])$ norm, as in \cite[eq. (5.52)]{ChK}.
Formulas \eqref{flaphi} give  the asymptotics in \eqref{psias} under the same normalization. 

Asymptotics \eqref{Dima} is obtained by integrating \eqref{phifla}, which gives 
$$
\int_0^1 \varphi_n(x)dx = C \nu_n^{-1}\big(1+O(\nu_n^{-1})\big),\quad n\to\infty,
$$
where $C\nu_n^{-1}$ is the integral of expression \eqref{phias} without the residual. 
Since this expression does not depend on  $\beta$, the constant factor $C$ must coincide with the value, which is obtained for $\beta=0$. 
In other words, the sequence of integrals $\int_0^1 \varphi_n(x)dx$ for the fOU process
and the fBm  have the same first order asymptotics. Hence the constant in \eqref{Dima} coincides with     
\cite[eq. (5.53)]{ChK}. 

\end{proof}

\subsubsection{Passing to the natural enumeration} 
The enumeration introduced in Lemma \ref{lem5.8} may not coincide with the  {\em natural} enumeration, which puts the eigenvalues into the decreasing order. Note that substitution of expression \eqref{nuas} into formula \eqref{lambdanu} gives the sequence  
$\lambda_n$, which decreases in the already chosen enumeration. Hence, starting from some index, these two enumerations may differ only by 
a constant shift. To determine this shift we can use the calibration procedure, based on the continuity of the spectrum with respect to the 
parameter $\alpha$ and the already known asymptotics \eqref{OUeig} for the standard OU process, corresponding to  
$\alpha=1$. This calibration is carried out exactly as for the fBm in \cite[Section 5.1.7]{ChK}, which shows 
that the formulas \eqref{nuas} and \eqref{phias}-\eqref{psias} must be shifted by 1: replacing $n$ by $n-1$, and 
$\alpha$ by $2-2H$, the equality \eqref{fBmeigfn} is obtained by Theorem \ref{thm-fOU}. 
%
%
%\begin{multline*} 
%\varphi_n(x) \propto    \sqrt{2}
%\sin \Big(\nu_n x + \frac {2H-1}8 \pi -\arcsin \frac{\ell_H}{\sqrt{1+\ell_H^2}}
%\Big)
%\\
%-
% \frac {\sqrt{2H+1}} \pi\int_0^\infty \rho_0(u)
% \Big(    -e^{-u\nu_n x}  \frac{u-\ell_H}{\sqrt{1+\ell_H^2}} +(-1)^{n}e^{-(1-x)u\nu_n} \Big)du + n^{-1}r_n(x),
%\end{multline*}
% 
%$$
%\int_0^1 e^{\beta x}\varphi_n(x)dx  \propto -
%  \sqrt{\frac{2H+1}{1+\ell_H^2}} \nu_n^{-1}  \\
%\quad \text{and}\quad 
%\varphi_n(1) \propto  - (-1)^{n}  
%\sqrt{2H+1}  \big(1+O(n^{-1})\big).
%$$
%

\subsection{The case $H<\frac 1 2$} 
In this case the covariance function is given by the formula \eqref{OUKst} 
and the spectral problem has the form
$$
\int_0^1 \left(
\int_0^x e^{\beta(x-u)}\frac d{du}\int_0^y e^{\beta(y-v)}C_\alpha|u-v|^{1-\alpha}\sign(u-v) dvdu
\right)\varphi(y)dy = \lambda \varphi(x),
$$
where $C_\alpha:=1-\frac \alpha 2$.  Differentiating twice we get  
$$
\int_0^1 \left(
  \frac d{dx}\int_0^y e^{\beta(y-v)}C_\alpha|x-v|^{1-\alpha}\sign(x-v) dv \right)\varphi(y)dy
  +
\beta \lambda \varphi(x) = \lambda \varphi'(x), 
$$
which can be written as  
\begin{align*}
&
-\frac d{dx}\int_0^1 \left(
  \int_0^y e^{-\beta v}C_\alpha|x-v|^{1-\alpha}\sign(x-v) dv \right)
\frac{d}{dy}  \int_y^1 e^{\beta r}\varphi(r)dr  dy
\\
&
  +
\beta \lambda \varphi(x) = \lambda \varphi'(x).
\end{align*}
Integrating by parts gives  
$$
 \frac d{dx}\int_0^1  C_\alpha|x-y|^{1-\alpha}\sign(x-y) \psi(y) dy + \beta \lambda\varphi(x) = \lambda \varphi'(x),
$$
where $\psi(x)$ is defined as in \eqref{psidef}. Using the identity \eqref{psiphi}
we get the generalized spectral problem, cf. \eqref{geneig}, 
$$
\begin{aligned}
\frac d{dx}\int_0^1  &  C_\alpha|x-y|^{1-\alpha}\sign(x-y) \psi(y) dy    = 
 \lambda\Big(\beta^2 \psi(x)-\psi''(x)\Big) 
 , \quad x\in [0,1], \\
&
\psi(1) =  0, \; 
\psi'(0)    +\beta \psi(0)  =0.
\end{aligned}
$$
The rest of the proof is carried out as in the case  $H>\frac 12$. 

\section{Concluding remarks}
We obtained the exact asymptotics of the mean squared error in the estimation problem of the fractional OU process 
observed in the white noise of vanishing intensity $\eps\to 0$. Due to the scaling property \eqref{Kbeta}  the results remain
valid on the arbitrary {\em finite} time interval $[0,T]$,  and the leading asymptotic term does not depend on the interval length $T$.
Another interesting problem would be to find the limit of the estimation error as $T\to\infty$ with the noise intensity $\eps>0$ being fixed.    
This large time asymptotic analysis requires the spectral estimates of Theorem \ref{thm-fOU} to be uniform in $T$. 
Such uniformity does not follow from the proof and the large time problem would need a different approach.

%\bibliographystyle{unsrt}
%\bibliography{/Users/Pavel/Dropbox/Pasha_Masha/bibliography/fBm}

\begin{thebibliography}{10}

\bibitem{M08}
Yuliya~S. Mishura.
\newblock {\em Stochastic calculus for fractional {B}rownian motion and related
  processes}, volume 1929 of {\em Lecture Notes in Mathematics}.
\newblock Springer-Verlag, Berlin, 2008.

\bibitem{PT17}
Vladas Pipiras and Murad~S. Taqqu.
\newblock {\em Long-range dependence and self-similarity}.
\newblock Cambridge Series in Statistical and Probabilistic Mathematics, [45].
  Cambridge University Press, Cambridge, 2017.

\bibitem{CKM03}
Patrick Cheridito, Hideyuki Kawaguchi, and Makoto Maejima.
\newblock Fractional {O}rnstein-{U}hlenbeck processes.
\newblock {\em Electron. J. Probab.}, 8:no. 3, 14, 2003.

\bibitem{LS1}
Robert~S. Liptser and Albert~N. Shiryaev.
\newblock {\em Statistics of random processes. {I}}, volume~5 of {\em
  Applications of Mathematics (New York)}.
\newblock Springer-Verlag, Berlin, expanded edition, 2001.
\newblock General theory, Translated from the 1974 Russian original by A. B.
  Aries, Stochastic Modelling and Applied Probability.

\bibitem{KB61}
Rudolph Kalman and Richard Bucy.
\newblock New results in linear filtering and prediction theory.
\newblock {\em Journal of Basic Engineering}, 83(1):95--108, 1961.

\bibitem{LS2}
Robert~S. Liptser and Albert~N. Shiryaev.
\newblock {\em Statistics of random processes. {II}}, volume~6 of {\em
  Applications of Mathematics (New York)}.
\newblock Springer-Verlag, Berlin, expanded edition, 2001.
\newblock Applications, Translated from the 1974 Russian original by A. B.
  Aries, Stochastic Modelling and Applied Probability.

\bibitem{IKh}
Il{\cprime}dar~A. Ibragimov and Rafail~Z. Has'minskii.
\newblock {\em Statistical Estimation. Asymptotic Theory.}, volume~16 of {\em
  Applications of Mathematics}.
\newblock Springer, 1981.

\bibitem{Tsybakov}
Alexandre~B. Tsybakov.
\newblock {\em Introduction to nonparametric estimation}.
\newblock Springer Series in Statistics. Springer, New York, 2009.
\newblock Revised and extended from the 2004 French original, Translated by
  Vladimir Zaiats.

\bibitem{Ukai}
Seiji Ukai.
\newblock Asymptotic distribution of eigenvalues of the kernel in the
  {K}irkwood-{R}iseman integral equation.
\newblock {\em J. Mathematical Phys.}, 12:83--92, 1971.

\bibitem{ChK}
Pavel Chigansky and Marina Kleptsyna.
\newblock Exact asymptotics in eigenproblems for fractional {B}rownian
  covariance operators.
\newblock {\em Stochastic Process. Appl.}, 128(6):2007--2059, 2018.

\bibitem{Gahov}
Fyodor~D. Gakhov.
\newblock {\em Boundary value problems}.
\newblock Dover Publications, Inc., New York, 1990.
\newblock Translated from the Russian, Reprint of the 1966 translation.

\end{thebibliography}

\def\cprime{$'$} \def\cprime{$'$} \def\cydot{\leavevmode\raise.4ex\hbox{.}}
  \def\cprime{$'$} \def\cprime{$'$} \def\cprime{$'$}

\end{document}